\documentclass[11pt,twoside]{article}
\usepackage{amssymb}

\evensidemargin\oddsidemargin
\newtheorem{definition}{Definition}[section]

\newtheorem{theorem}[definition]{Theorem}
\newtheorem{proposition}[definition]{Proposition}

\newtheorem{remark}[definition]{Remark}

\newtheorem{examples}[definition]{Examples}
\def\gpd{\rightrightarrows}
\def\gr{G\gpd M}
\def\grr{G\times \R\gpd M\times \R}
\def\s{\alpha}
\def\t{\beta}
\def\mult{\sigma :G\to \R}
\def\Jacobi{\Lambda , E}
\def\sostJ{\# _{(\Lambda ,E)}}
\def\sostP{\#_\Lambda}
\def\QED{\hskip0.1em\hfill\null\ \null\nobreak\hfill
\kern3pt\lower1.8pt\vbox{\hrule\hbox
{\vrule\kern1pt\vbox{\kern1.7pt \hbox{$\scriptstyle
QED$}\kern0.2pt}\kern1pt\vrule}\hrule}}
\let\izq=\overleftarrow
\let\der=\overrightarrow
\def\exps{e^{\sigma (g)}}

\def\inv{^{-1}}

\font\ddpp=msbm9  scaled \magstep 1  
\newfont{\msi}{msbm8 scaled \magstephalf}  

\def\Rp{\hbox{\msi R}}           
\def\R{\hbox{\ddpp R}}               

\def\lcf{\lbrack\! \lbrack}
\def\rcf{\rbrack\! \rbrack}
\newcommand\prueba {\mbox{{\em Proof: }}}

\begin{document}

\label{xyzt} \leavevmode\vadjust{\vskip -50pt}


%
%

\pagestyle{myheadings} \markboth {{\small\sc l.c.s. groupoids}}
{{\small\sc D. Iglesias-Ponte, J.C. Marrero}}

%
%

\thispagestyle{empty}

\begin{center}
{\Large\bf Locally conformal symplectic groupoids } \vskip 10mm

%
%

{\large\bf David Iglesias-Ponte and J.C. Marrero } \vskip 5mm {\it
Departamento de Matem\'atica Fundamental, Universidad de La Laguna

\smallskip

emails: diglesia@ull.es, jcmarrer@ull.es }
\end{center}

\bigskip

%
%

\begin{abstract}
\parindent0pt\noindent
Locally conformal symplectic (l.c.s.) groupoids are introduced as
a generalization of sym\-plec\-tic group\-oids. We obtain some
examples and we prove that l.c.s. groupoids are examples of Jacobi
groupoids in the sense of \cite{IM}. Finally, we describe the Lie
algebroid of a l.c.s. groupoid.



\end{abstract}

%
%
\section{Introduction}
A symplectic groupoid is a Lie groupoid $\gr$ endowed with a
symplectic structure $\Omega$ for which the graph of the partial
multiplication is a lagrangian submanifold in the symplectic
manifold $(G\times G\times G,\Omega \oplus \Omega \oplus -\Omega
)$. This interesting class of groupoids, which was introduced in
\cite{CDW}, arises in the integration of arbitrary Poisson
manifolds. In fact, if $(\gr ,\Omega )$ is a symplectic groupoid
then there exists a Poisson structure $\Lambda _0$ on $M$ and the
Lie algebroid $AG$ is isomorphic to the cotangent Lie algebroid
$T^\ast M$.

An interesting generalization of symplectic groupoids, as well as
of Drinfeld's Poisson-Lie groups, are Poisson groupoids \cite{We}.

On the other hand, a non-degenerate 2-form on a manifold $M$ of
even dimension is said to be locally conformal symplectic (l.c.s.)
if it is conformally related with a symplectic structure in some
neighbourhood of every point of $M$ (see \cite{GL,Vai1}). L.c.s.
manifolds are interesting examples of Jacobi manifolds \cite{GL}
and, in addition, they play an important role in the study of some
dynamical systems, particularly conformally Hamiltonian systems
(see \cite{LW}). The aim of this paper is to introduce the notion
of a l.c.s. groupoid and to study some of its properties.

The paper is organized as follows. In Section 2, we recall the
definition of a l.c.s. manifold and its relation with Jacobi
manifolds. In Section 3, we prove that if $G$ is a contact
groupoid over a manifold $M$ then it is possible to introduce a
l.c.s. structure and a Lie groupoid structure on $G\times \R$
(whose space of identity elements is $M\times \R$) and these
structures are compatible in a certain way. These results motivate
the definition of a l.c.s. groupoid in Section 4. Symplectic
groupoids are l.c.s. groupoids and, furthermore, we prove, in
Section 4, that a l.c.s. groupoid is a particular example of the
so-called Jacobi groupoids, which were first introduced in
\cite{IM} as a generalization of Poisson groupoids. Finally, in
Section 5, we describe the Lie algebroid of a l.c.s. groupoid.

In this paper, we will use the definitions, notation and
conventions introduced in our previous paper \cite{IM} (see
Sections 1 and 2 in \cite{IM}).
\section{Locally conformal symplectic manifolds}
A manifold $M$ is said to be {\it locally (globally) conformal
symplectic (l.(g.)c.s.) manifold} if it admits a non-degenerate
$2$-form $\Omega$ and a closed (exact) 1-form $\omega$ such that
\begin{equation}\label{0'}
\delta \Omega = \omega \wedge \Omega.
\end{equation}
The $1$-form $\omega$ is called the {\it Lee 1-form} of $M$. It is
obvious that the l.c.s. manifolds with Lee 1-form identically zero
are just the symplectic manifolds (see, for example,
\cite{GL,Vai1}).

L.c.s manifolds are examples of Jacobi manifolds, i.e., if
$(M,\Omega ,\omega)$ is a l.c.s. mani\-fold then there exists a
2-vector $\Lambda$ and a vector field $E$ on $M$ such that
$[\Lambda ,\Lambda ]=2E\wedge \Lambda$ and $[E,\Lambda ]=0$. In
fact, the Jacobi structure $(\Lambda,E)$ is given by
\begin{equation}\label{lcs->jacobi}
\Lambda(\alpha,\beta) =
\Omega(\flat^{-1}(\alpha),\flat^{-1}(\beta)), \makebox[1cm]{} E =
\flat^{-1}({\omega}) \; ,
\end{equation}
 for all $\alpha , \beta \in
\Omega^{1}(M)$, where $\flat : \mathfrak{X}(M) \longrightarrow
\Omega^1(M)$ is the isomorphism of $C^{\infty}(M,\R)$-modules
defined by $\flat (X) = i(X)\Omega$ (see \cite{GL}).

\section{An example}
Let $(\gr ,\eta ,\sigma )$ be a {\em contact groupoid} over a
manifold $M$, that is, $\gr$ is a Lie groupoid over $M$ with
structural functions $\s$, $\t$, $m$ and $\epsilon$, $\eta \in
\Omega ^1(G)$ is a contact 1-form on $G$ and $\sigma :G\to \R$ is
an arbitrary function such that if $\oplus _{TG}$ is the partial
multiplication in the tangent Lie groupoid $TG\gpd TM$, then
\begin{equation}\label{contgroup}
\eta _{(gh)}(X_g \oplus _{TG} Y_h )=\eta _{(g)}(X_g)+\exps \eta
_{(h)}(Y_h),
\end{equation}
for $(g,h)\in G ^{(2)}$ and $(X_g,Y_h)\in T_{(g,h)}G^{(2)}$ (see
\cite{D2,KS}).

If $(\gr ,\eta ,\sigma )$ is a contact groupoid then, using the
associativity of $\oplus _{TG}$, we deduce that $\mult$ is a
multiplicative function, that is, $\sigma (gh)=\sigma (g)+\sigma
(h)$, for $(g,h)\in G^{(2)}$. This implies that
\begin{equation}\label{2'}
\sigma \circ \epsilon =0.
\end{equation}
Using (\ref{contgroup}), we also deduce that (see \cite{D2})
\begin{equation}\label{utiles}
\xi (\sigma )=0,
\end{equation}
\begin{equation}\label{2.2'}
\begin{array}{lll} (\delta \eta )_{gh} (X_g\oplus
_{TG}Y_h,X_g'\oplus _{TG}Y'_h)&\kern-8pt=&\kern-8pt (\delta \eta)
_g (X_g,X_g')+e^{\sigma (g)}(\delta \eta )_h(Y_h,Y'_h)
\\&&\kern-20pt+e^{\sigma (g)}(X_g(\sigma )\eta _h(Y'_h)- X'_g(\sigma )\eta
_h(Y_h)),\end{array}
\end{equation}
for $(X_g,Y_h),(X'_g,Y'_h)\in T_{(g,h)}G^{(2)}$, $\xi$ being the
Reeb vector field of the contact structure $\eta$.

Now, we will obtain a l.c.s structure and a Lie groupoid structure
on the manifold $G\times\R$, both structures compatible in a
certain way.

First of all, take the Lie groupoid $\gr$ and the multiplicative
function $\mult$. Then, using the multiplicative character of
$\sigma$, we can define a right action of $\gr$ on the canonical
projection $\pi _1:M\times \R \to M$ as follows
\begin{equation}\label{accion}
(x,t)\cdot g=(\s (g),\sigma (g)+t)
\end{equation}
for $(x,t)\in M\times \R$ and $g\in G$ such that $\t(g)=x.$ Thus,
we have the corresponding action groupoid $(M\times \R)\ast G\gpd
M\times \R$ over $M\times\R$, where
\[
(M\times \R)\ast G =\{ ((x,t),g)\in (M\times\R )\times G\, / \,
\beta (g)=x\}
\]
(see \cite{HM} for the general definition of an action Lie
groupoid; see also \cite{IM}). Moreover, it is not difficult to
prove that $(M\times \R)\ast G$ may be identified with the product
manifold $G\times \R$ and, under this identification, the
structural functions of the Lie groupoid are given by
\begin{equation}\label{GxR}
\begin{array}{rcll}
\s_\sigma(g,t)&=&(\s(g), \sigma(g) + t),&\mbox{ for } (g,t)\in
G\times \R,\\ \t_\sigma(h,s)&=&(\t(h),s),& \mbox{ for } (h,s)\in
G\times \R,\\ m_\sigma((g,t),(h,s))&=&(gh,t),&\mbox{ if } \s
_\sigma(g,t)=\t _\sigma(h,s),\\ \epsilon_\sigma(x,t)&=&(\epsilon
(x), t),& \mbox{ for } (x,t)\in M\times \R .
\end{array}
\end{equation}
From (\ref{GxR}) and the definition of the tangent groupoid (see,
for instance, \cite{IM}), it follows that the projections
$(\s_\sigma)^T,(\t_\sigma)^T$, the inclusion $(\epsilon_\sigma)^T$
and the partial multiplication $\oplus_{T(G\times \Rp)}$ of the
tangent groupoid $T(G\times \R)\gpd T(M\times \R)$ are given by
\begin{equation}\label{T(GxR)}
\begin{array}{l}
\kern-10pt(\s_\sigma)^T(X_g + \lambda\frac{\partial}{\partial
t}_{|t}) =\s^T(X_g) + (\lambda + X_g(\sigma))\frac{\partial
}{\partial t}_{|t+\sigma(g)},\\ \kern-10pt(\t_\sigma)^T(Y_h +
\mu\frac{\partial}{\partial t}_{|s})= \t^T(Y_h) + \mu
\frac{\partial}{\partial t}_{|s},
\\ \kern-10pt (X_g +
\lambda\frac{\partial}{\partial t}_{|t})\oplus_{T(G\times \Rp)}(
Y_h + \mu\frac{\partial}{\partial t}_{|s})= X_g\oplus_{TG}Y_h +
\lambda\frac{\partial }{\partial t}_{|t}, \\ \kern-10pt
 (\epsilon_\sigma)^T(X_x
+ \lambda\frac{\partial}{\partial t}_{|t})=\epsilon^T(X_x) +
\lambda\frac{\partial}{\partial t}_{|t},
\end{array}
\end{equation}
for $X_g + \lambda\frac{\partial}{\partial t}_{|t}\in
T_{(g,t)}(G\times \R)$, $Y_h + \mu\frac{\partial}{\partial
t}_{|s}\in T_{(h,s)}(G\times \R)$ and $X_x + \lambda
\frac{\partial }{\partial t}_{|t} \in T_{(x,t)}(M\times \R)$ (see
Section 5.3.1 in \cite{IM}).

On the other hand, if $A(G\times\R )$ is the Lie algebroid of the
Lie groupoid $G\times\R\gpd M\times\R$, then using (\ref{GxR}) and
the definition of the cotangent group\-oid (see \cite{CDW,IM}), we
deduce that the projections
$\widetilde{\s_\sigma},\widetilde{\t_\sigma}$, the inclusion
$\widetilde{\epsilon_\sigma}$ and the partial multiplication
$\oplus_{T^*(G\times \Rp)}$ in the cotangent groupoid $T^\ast
(G\times \R)\gpd A^\ast (G\times \R)$ are defined by
\begin{equation}\label{T*(GxR)}
\begin{array}{l}
\widetilde{\s_\sigma}(\mu _g + \gamma \delta
t_{|t})=(\tilde{\s}(\mu _g),\sigma(g) + t),\\[2pt]
\widetilde{\t_\sigma}(\nu_h + \zeta \delta
t_{|s})=(\tilde{\t}(\nu_h)-\zeta (\delta\sigma)_{\epsilon (\t(g))}
,s),\\[2pt] (\mu_g +\gamma \delta t_{|t})\oplus_{T^\ast(G\times
\Rp)}( \nu_h + \zeta \delta t_{|s})= (\mu_g +
\zeta(\delta\sigma)_g)\oplus_{T^\ast G}\nu_h + (\gamma +
\zeta)\delta t_{|t}\\[2pt]
\widetilde{\epsilon_\sigma}(\mu_x,t)=\tilde\epsilon(\mu_x) +
0\,\delta t_{|t},
\end{array}
\end{equation}
for $\mu_g + \gamma \delta t_{|t}\in T_{(g,t)}^\ast(G\times \R)$,
$\nu_h + \zeta \delta t_{|s}\in T_{(h,s)}^\ast(G\times \R)$ and
$(\mu_x,t)\in A^\ast _{(x,t)}(G\times \R)\cong A_x^*G\times \R$,
where $\tilde{\s}$, $\tilde{\t}$, $\oplus _{T^\ast G}$ and
$\tilde{\epsilon}$ are the structural functions of the cotangent
Lie groupoid $T^\ast G\gpd A^\ast G$ (see Section 5.3.1 in
\cite{IM}).

Next, we define on $G\times \R$ the 2-form $\Omega$ and the 1-form
$\omega$ given by
\begin{equation}\label{lcs-can}
\Omega = - (\bar{\pi}^\ast _1 (\delta \eta )+\bar{\pi}^\ast _2
(\delta t)\wedge \bar{\pi}^\ast _1(\eta )),\quad \omega
=-\bar{\pi}^\ast _2 (\delta t),
\end{equation}
where $\bar{\pi}_1:G\times \R\to G$ and $\bar{\pi}_2:G\times \R\to
\R$ are the canonical projections. Thus, we have that $(G\times
\R, \Omega ,\omega )$ is a locally conformal symplectic manifold
of the first kind in the sense of \cite{Vai1}. Note that if
$(\Lambda ,E)$ is the Jacobi structure on $G\times\R$ associated
with the l.c.s. structure $(\Omega ,\omega )$ then
\begin{equation}\label{E=-Reeb}
E=-\xi.
\end{equation}
Now, using (\ref{contgroup})-(\ref{2.2'}) and
(\ref{GxR})-(\ref{E=-Reeb}), we prove the following result.
\begin{proposition}\label{prop-lcs}
Let $(\gr ,\eta ,\sigma )$ be a contact groupoid. If $\grr$ is the
Lie groupoid with structural functions given by (\ref{GxR}), the
pair $(\Omega ,\omega )$ is defined by (\ref{lcs-can}),
$\bar{\sigma}= \sigma \circ \bar{\pi}_1$ is the pull-back of the
multiplicative function $\sigma$ by the canonical projection
$\bar{\pi}_1:G\times\R\to G$ and $\theta$ is the 1-form on
$G\times\R$ defined by $\theta =e^{\bar{\sigma}} (\delta
\bar{\sigma}-\omega )$, then we have:
\begin{itemize}
\item[i)] $m_\sigma ^\ast \Omega =\pi _1^\ast \Omega +
e^{(\bar{\sigma}\circ \pi _1)}\pi _2^\ast \Omega$;
\item[ii)] $\widetilde{\s_\sigma}\circ \omega =0$,
$\widetilde{\t_\sigma}\circ \theta =0$;
\item[iii)] $m_\sigma ^\ast \omega =\pi _1^\ast \omega$,
$m_\sigma^\ast \theta= e^{(\bar{\sigma}\circ \pi _1)}\pi
_2^\ast\theta$;
\item[iv)] $\Lambda (\omega ,\theta )=0$, $(\theta + \omega
-\widetilde{\epsilon _\sigma }\circ \widetilde{\beta _\sigma}\circ
\omega )\circ \epsilon =0$;
\end{itemize}
where $\pi _i :G^{(2)}\to G$, $i=1,2$, are the canonical
projections and $(\Lambda ,E)$ is the Jacobi structure on $G\times
\R$ associated with the l.c.s. structure $(\Omega ,\omega )$.
\end{proposition}
\section{\kern-7pt Locally conformal symplectic
\kern-.4pt and \kern-.4pt Jacobi groupoids} Motivated by
Proposition \ref{prop-lcs}, we introduce the following definition.
\begin{definition}\label{def-lcs}
Let $\gr$ be a Lie groupoid with structural functions $\s$, $\t$,
$m$ and $\epsilon$, $(\Omega ,\omega )$ be a l.c.s. structure on
$G$, $\mult$ be a multiplicative function and $\theta$ be the
1-form on $G$ defined by
\begin{equation}\label{theta-def}
\theta =e^{\sigma} (\delta \sigma -\omega ).
\end{equation}
Then, $(\gr,\Omega, \omega ,\sigma )$ is a {\em l.c.s. groupoid}
if the following properties hold:
\begin{equation}\label{lagrangiano}
m^\ast \Omega =\pi _1^\ast \Omega + e^{(\sigma \circ \pi _1)}\pi
_2^\ast \Omega;
\end{equation}
\begin{equation}\label{14'}
\tilde{\alpha}\circ \omega =0, \quad \tilde{\beta}\circ \theta =0;
\end{equation}
\begin{equation}\label{14''}
m^\ast \omega =\pi _1^\ast \omega,\quad m^\ast \theta= e^{(\sigma
\circ \pi _1)}\pi _2^\ast\theta;
\end{equation}
\begin{equation}\label{14'''}
\Lambda (\omega ,\theta )=0,\quad (\theta + \omega
-\tilde{\epsilon}\circ \tilde{\beta}\circ \omega )\circ \epsilon
=0;
\end{equation}
where $\pi _i :G^{(2)}\to G$, $i=1,2$, are the canonical
projections, $(\Lambda ,E)$ is the Jacobi structure associated
with the l.c.s. structure $(\Omega ,\omega )$ and $\tilde{\s}$,
$\tilde{\t}$, $\oplus _{T^\ast G}$ and $\tilde{\epsilon}$ are the
structural functions of the cotangent groupoid $T^\ast G\gpd
A^\ast G$.
\end{definition}
\begin{examples}\label{ejemplos}
{\rm {\it i)} If $(\gr ,\eta ,\sigma )$ is a contact groupoid
then, by Proposition \ref{prop-lcs}, $(\grr, \Omega ,\omega ,
\bar{\sigma})$ is a l.c.s. groupoid, where the structural
functions of $\grr$ are defined by (\ref{GxR}), the pair $(\Omega
,\omega )$ is given by (\ref{lcs-can}) and the multiplicative
function $\bar{\sigma}$ is $\bar{\sigma}= \sigma \circ
\bar{\pi}_1$, $\bar{\pi}_1:G\times\R\to G$ being the canonical
projection.

{\it ii)} A Lie groupoid $\gr$ is said to be {\em symplectic} if
$\,G$ admits a symplectic 2-form $\Omega$ in such a way that the
graph of the partial multiplication in $G$ is a Lagrangian
submanifold of the symplectic manifold $(G\times G\times G,\Omega
\oplus \Omega \oplus (-\Omega ))$ (see \cite{CDW}). This is
equivalent to say that $\Omega$ satisfies the condition $ m^\ast
\Omega =\pi _1^\ast \Omega + \pi _2^\ast \Omega.$ Therefore, we
conclude that $(\gr ,\Omega )$ is a symplectic groupoid if and
only if $(\gr, \Omega ,0 ,0)$ is a l.c.s. groupoid. }
\end{examples}
Next, we will give the relation between l.c.s. groupoids and
Jacobi groupoids. First, we will recall the definition of a Jacobi
groupoid.

Let $\gr$ be a Lie groupoid and $\mult$ be a multiplicative
function. Then, $TG\times \R$ is a Lie groupoid over $TM\times\R$
with structural functions given by (see Section 3 in \cite{IM})
\begin{equation}\label{TGR}
\begin{array}{l} (\s ^T)_\sigma (X_g,\lambda )=(\s ^T(X_g),
X_g(\sigma )+\lambda ),\mbox{ for }(X_g,\lambda )\in T_gG\times\R,
\\(\t ^T)_\sigma (Y_h,\mu )=(\t ^T(Y_h),\mu ),
\mbox{ for }(Y_h,\mu )\in T_hG\times\R ,\\ (X_g,\lambda )\oplus
_{TG\times \Rp}(Y_h,\mu )=(X_g\oplus_{TG} Y_h,\lambda ), \mbox{ if
}(\s ^T)_\sigma (X_g,\lambda )=(\t ^T)_\sigma (Y_h,\mu ),\\
(\epsilon ^T)_\sigma (X_x,\lambda )=(\epsilon ^T(X_x),\lambda
),\mbox{ for }(X_x,\lambda )\in T_xM\times\R .
\end{array}
\end{equation}
On the other hand, if $AG$ is the Lie algebroid of $G$ then
$T^\ast G\times \R$ is a Lie groupoid over $A^\ast G$ with
structural functions defined by (see Section 3 in \cite{IM})
\begin{equation}\label{T*GR}
\begin{array}{l}
\tilde{\s}_\sigma(\mu _g,\gamma )=e^{-\sigma(g)} \tilde{\s}(\mu
_g), \mbox{ for }(\mu _g,\gamma )\in T_g^\ast G\times \R,\\
\tilde{\t}_\sigma(\nu _h,\zeta )=\tilde{\t}(\nu _h)-\zeta \,
(\delta \sigma )_{\epsilon (\t (h))}{}_{|A_{\t (h)}G} ,\mbox{ for
}(\nu _h,\zeta )\in T_h^\ast G\times \R,\\ (\mu _g,\gamma )\oplus
_{T^\ast G\times \Rp}(\nu _h,\zeta ) =\Big ( (\mu _g +e^{\sigma
(g)}\zeta \,(\delta \sigma ) _g)\oplus _{T^\ast G}(e^{\sigma
(g)}\nu _h),\gamma +e^{\sigma (g)}\zeta \Big )\\
\tilde{\epsilon}_\sigma (\mu _x)=(\tilde{\epsilon}(\mu
_x),0),\mbox{ for }\mu _x\in A^\ast _xG .
\end{array}
\end{equation}
\begin{definition}\cite{IM}\label{def-Jacgroup}
Let $\gr$ be a Lie groupoid, $(\Jacobi )$ be a Jacobi structure on
$G$ and $\sigma :G\to \R$ be a multiplicative function. Then,
$(\gr,\Jacobi ,\sigma )$ is a {\em Jacobi groupoid} if the
homomorphism $\sostJ:T^\ast G\times \R\to TG\times \R$ given by
$$\sostJ (\mu _g,\gamma )=(\sostP(\mu _g)+\gamma \,E_g,-\mu
_g(E_g))$$is a morphism of Lie groupoids over some map $\varphi
_0:A^\ast G\to TM\times \R$, where the structural functions of the
Lie groupoid structure on $T^\ast G\times \R\gpd A^\ast G$
(respectively, $TG\times \R \gpd TM\times \R$) are given by
(\ref{T*GR}) (respectively, (\ref{TGR})).
\end{definition}
\begin{remark}
{\rm A Poisson groupoid is a Jacobi groupoid $(\gr ,\Jacobi
,\sigma )$ with $E=0$ and $\sigma =0$ (see \cite{IM}).}
\end{remark}

A characterization of a Jacobi groupoid is the following one. If
$\gr$ is a Lie groupoid and $\mult$ is a multiplicative function
then $T^\ast G$ is a Lie groupoid over $A^\ast G$ with structural
functions given by
\begin{equation}\label{T*Gsigma}
\begin{array}{l}
\tilde{\s}^\ast _\sigma(\mu _g)=e^{-\sigma(g)} \tilde{\s}(\mu _g),
\mbox{ for }\mu _g\in T_g^\ast G,\\ \tilde{\t}^\ast _\sigma(\nu
_h)=\tilde{\t}(\nu _h),\mbox{ for }\nu _h\in T_h^\ast G,\\ \mu
_g\oplus ^\sigma _{T^\ast G}\nu _h =\mu _g \oplus _{T^\ast
G}e^{\sigma (g)}\nu _h, \mbox{ if }\tilde{\s}^\ast _\sigma(\mu
_g)=\tilde{\t}^\ast _\sigma(\nu _h),\\ \tilde{\epsilon}^\ast
_\sigma (\mu _x)=\tilde{\epsilon} (\mu _x), \mbox{ for }\mu _x\in
A^\ast _xG.
\end{array}
\end{equation}
We call this Lie groupoid the $\sigma${\em -cotangent groupoid}.
Note that the canonical inclusion $T^\ast G\to T^\ast G\times\R$,
$\mu _g\mapsto (\mu _g,0)$, is a monomorphism of Lie groupoids.
\begin{proposition}\label{charac-Jacgroup}
Let $\gr$ be a Lie groupoid, $(\Jacobi )$ be a Jacobi structure on
$G$ and $\sigma :G\to \R$ be a multiplicative function. Then,
$(\gr,\Jacobi ,\sigma )$ is a {\em Jacobi groupoid} if and only if
the following conditions hold:
\begin{itemize}
\item[{\it i)}] $\sostP :T^\ast G\to TG$ is a Lie groupoid morphism
over some map $\tilde{\varphi} _0:A^\ast G\to TM$ from the
$\sigma$-cotangent groupoid $T^\ast G\gpd A^\ast G$ to the tangent
Lie groupoid $TG\gpd TM$.
\item[{\it ii)}] $E$ is a right-invariant vector field on $G$ and
$E(\sigma )=0$.
\item[{\it iii)}]  If $X_0\in \Gamma (AG)$
is the section of the Lie algebroid $AG$ satisfying
$E=-\der{X_0}$, we have that
\[
\sostP (\delta \sigma )=\der{X_0}-e^{-\sigma}\izq{X_0}.
\]
\end{itemize}
\end{proposition}
\prueba Suppose that $(\gr ,\Jacobi ,\sigma )$ is a Jacobi
groupoid. Then, proceeding as in the proof of Proposition 4.4 in
\cite{IM}, we deduce that {\it i), ii)} and {\it iii)} hold.

A similar computation proves the converse.\QED

Now, we will show that a l.c.s. symplectic groupoid is a
particular example of a Jacobi groupoid.
\begin{theorem}\label{lcs->Jac}
Let $\gr$ be a Lie groupoid, $(\Omega ,\omega)$ be a l.c.s.
structure on $G$ and  $\mult$ be a multiplicative function. If
$(\Jacobi )$ is the Jacobi structure associated with the l.c.s.
structure $(\Omega ,\omega )$ then $(\gr ,\Omega ,\omega ,\sigma
)$ is a l.c.s. groupoid if and only if $(\gr ,\Jacobi ,\sigma )$
is a Jacobi groupoid.
\end{theorem}
\prueba Assume that $(\gr ,\Omega ,\omega ,\sigma )$ is a l.c.s.
groupoid. If $\mu _g\in T^\ast _gG$ and $\nu _h\in T^\ast _hG$
satisfy the relation $\tilde{\s}^\ast _\sigma (\mu
_g)=\tilde{\t}^\ast _\sigma (\nu _h)$, then using
(\ref{lcs->jacobi}), (\ref{lagrangiano}) and the definition of the
partial multiplication $\oplus _{TG}$ in the tangent Lie groupoid
$TG\gpd TM$, we obtain that $\s ^T (\sostP (\mu _g))=\t ^T(\sostP
(\nu _h))$ and, in addition,
\[
\begin{array}{l}
\kern-70pt(i(\sostP (\mu _g\oplus ^\sigma _{T^\ast G} \nu
_h))\Omega _{(gh)})(X_g\oplus _{TG}Y_h)\\[4pt]=\, \, (i(\sostP
(\mu _g)\oplus _{TG}\sostP (\nu _h) )\Omega _{(gh)})(X_g\oplus
_{TG}Y_h),
\end{array}
\]
for $(X_g, Y_h)\in T_{(g,h)}G^{(2)}$. Thus (see
(\ref{lcs->jacobi})), it follows that $\sostP (\mu _g\oplus
^\sigma _{T^\ast G} \nu _h)=\sostP (\mu _g)\oplus _{TG} \sostP
(\nu _h)$ and, therefore, the map $\sostP:T^\ast G\to TG$ is a Lie
groupoid morphism over some map $\tilde{\varphi}_0:A^\ast G\to
TM$, between the $\sigma$-cotangent groupoid $T^\ast G\gpd A^\ast
G$ and the tangent groupoid $TG\gpd TM$. In particular, this
implies that
\begin{equation}\label{relaciones}
\s ^T\circ \sostP=\tilde{\varphi} _0\circ \tilde{\s}_\sigma
^\ast,\qquad  \t ^T\circ \sostP=\tilde{\varphi} _0\circ
\tilde{\t}_\sigma ^\ast.
\end{equation}
Now, from (\ref{lcs->jacobi}), we deduce that $E=-\sostP(\omega
)$. Using this relation, (\ref{14'}) and (\ref{relaciones}), we
have that the vector field $E$ is $\s$-vertical

Next, suppose that $(g,h)\in G^{(2)}$ and denote by $R_h:G_{\t
(h)} \to G_{\s (h)}$ the right-translation by $h$. Then,
(\ref{lcs->jacobi}) and (\ref{14''}) imply that
\[
(i(E_{(gh)})\Omega _{(gh)})(X_g\oplus _{TG} Y_h)= (i((R_h)_\ast ^g
(E_{(g)}))\Omega _{(gh)})(X_g\oplus _{TG} Y_h),
\]
for $(X_g, Y_h)\in T_{(g,h)}G^{(2)}$. Consequently, $E$ is a
right-invariant vector field and there exists $X_0\in \Gamma (AG)$
such that $E=-\der{X}_0$.

On the other hand, if $X_{e^{\sigma}}$ is the hamiltonian vector
field of the function $e^{\sigma}$,
$X_{e^{\sigma}}=e^{\sigma}\sostP (\delta \sigma )+e^{\sigma }E$,
it is clear that $X_{e^{\sigma}}=\sostP (\theta )$. Using this
equality, (\ref{14'}), (\ref{14''}), (\ref{relaciones}) and
proceeding as in the proof of the fact that $E$ is
right-invariant, we conclude that $X_{e^{\sigma}}$ is a
left-invariant vector field. Moreover, if $x$ is a point of $M$,
then relation (\ref{14'''}) implies that $X_{e^{\sigma}}(\epsilon
(x) )=-\izq{X}_0(\epsilon (x))$. Thus, $\sostP (\delta \sigma )=
\der{X}_0- e^{-\sigma}\izq{X}_0$.

Finally, since $\Lambda (\omega ,\theta )=0$ and $E=-\sostP
(\omega )$, we obtain that $E(\sigma)=0$. Therefore, $(\gr
,\Jacobi ,\sigma )$ is a Jacobi groupoid.

In a similar way, we prove the converse.\QED

\begin{remark}
{\rm Using Theorem \ref{lcs->Jac} we directly deduce that a
symplectic group\-oid is a Poisson groupoid. This result was
proved in \cite{We}. }
\end{remark}
\section{The Lie algebroid of a l.c.s.
groupoid} Let $(\gr ,\Omega ,\omega ,\sigma )$ be a l.c.s.
groupoid and $\theta$ the 1-form on $G$ given by
(\ref{theta-def}). Then, the 1-form $e^{-\sigma}\theta$ is closed
and since $\tilde{\t}\circ\theta =0$, it follows that $\theta$ is
basic with respect to the projection $\s$. Thus, there exists a
unique 1-form $\theta _0$ on $M$ such that $\s ^\ast \theta _0
=e^{-\sigma}\theta $. It is clear that $\theta _0$ is closed.

Now, denote by $(\Jacobi )$ the Jacobi structure on $G$ associated
with the l.c.s. structure $(\Omega ,\omega )$. Then,
$\sostP(\theta)$ is the hamiltonian vector field $X_{e^{\sigma}}$
of the function $e^\sigma$. Moreover, from Theorem \ref{lcs->Jac},
and using the results in \cite{IM} (see Proposition 5.6 in
\cite{IM}), we deduce that there exists a Jacobi structure
$(\Lambda _0,E_0)$ on $M$ in such a way that the couple $(\s
,e^\sigma)$ is a conformal Jacobi morphism between the Jacobi
manifolds $(G,\Jacobi)$ and $(M,\Lambda_0,E_0)$. This implies that
\[
\begin{array}{lll}
\# _{\Lambda _0}{}_{(\s (g))} (\theta _0 (\s
(g)))\kern-.25pt&=&\kern-.25pt e ^{\sigma (g)} (\s ^T _g \circ
\sostP {}_{(g)} \circ (\s ^T_g)^\ast ) (\theta _0 (\s(g))
)\kern-.25pt \\[7pt] &=& \kern-.25pt \s ^T_g
(X_{e^{\sigma}}(g))\kern-.25pt =\kern-.25pt E_0 (\s(g)),
\end{array}
\]
for $g\in G$, where $(\s ^T_g)^\ast :T^\ast _{\s (g)}M\to T^\ast
_gG$ is the adjoint map of the tangent map $\s ^T_g:T_gG\to T_{\s
(g)}M$. Therefore, we have proved the following result.
\begin{proposition}\label{enM}
Let $(\gr ,\Omega ,\omega ,\sigma )$ be a l.c.s. groupoid and
$\theta$ the 1-form on $G$ given by (\ref{theta-def}). Then, there
exists a unique 1-form $\theta _0$ on $M$ such that $\s ^\ast
\theta _0=e^{-\sigma}\theta$. Furthermore, $\theta _0$ is closed
and $\# _{\Lambda _0} (\theta _0)=E_0$.
\end{proposition}
Next, we will describe the Lie algebroid associated with a l.c.s.
groupoid.
\begin{theorem}\label{algebroide}
Let $(\gr ,\Omega ,\omega ,\sigma )$ be a l.c.s. groupoid, $AG$ be
the Lie algebroid of $G$, $(\Jacobi )$ be the Jacobi structure on
$G$ associated with the l.c.s. structure $(\Omega ,\omega )$ and
$(\Lambda _0,E_0)$ be the corresponding Jacobi structure on $M$.
Then, the map $\Psi :\Omega ^1(M)\to \mathfrak X _L(G)$ between
$\Omega ^1(M)$ and the space of left-invariant vector fields on
$G$ defined by $\Psi (\mu )=e^{\sigma}\sostP (\s ^\ast \mu )$
induces an isomorphism between the vector bundles $T^\ast M$ and
$AG$. Under this isomorphism, the Lie bracket on $\Gamma (AG)\cong
\mathfrak X _L(G)$ and the anchor map of $AG$ are given by
\[
\begin{array}{lll}
\lcf \mu ,\nu \rcf _{(\Lambda _0,E_0,\theta _0)}&=&{\cal
L}_{\#_{\Lambda _0}(\mu )}\nu\kern-2pt-\kern-2pt{\cal
L}_{\#_{\Lambda _0}(\nu )}\mu\kern-2pt -\kern-2pt \delta (\Lambda
_0(\mu ,\nu ))\\ & & -i(E_0)(\mu \kern-1pt\wedge\kern-1pt
\nu)-\Lambda _0(\mu ,\nu )\theta _0,\\ \widetilde{\#}_{(\Lambda
_0,E_0,\theta _0)}(\mu )&=&\#_{\Lambda _0}(\mu ),
\end{array}
\]
for $\mu,\nu\in \Omega ^1(M)$, where $\theta _0$ is the 1-form on
$M$ considered in Proposition \ref{enM}.
\end{theorem}
\prueba Let $\mu$ be a 1-form on $M$. Since the map $\sostP:T^\ast
G\to TG$ is a morphism between the $\sigma$-cotangent groupoid and
the tangent groupoid $TG\gpd TM$, we obtain the that vector field
$\tilde{X} =\Psi (\mu )$ is $\t$-vertical. Moreover, if $(g,h)\in
G^{(2)}$ and $L_g:G^{\s (g)}\to G^{\t (g)}$ is the
left-translation by $g$ then, using (\ref{lagrangiano}), we deduce
that
\[
(i(\tilde{X}_{(gh)})\Omega _{(gh)})(Y_g\oplus _{TG}
Z_h)=(i((L_g)_\ast ^h (\tilde{X}_{(h)}))\Omega _{(gh)})
(Y_g\oplus_{TG} Z_h).
\]
for $(Y_g, Z_h)\in T_{(g,h)}G^{(2)}$. This proves that
$\tilde{X}\in \mathfrak X _L(G)$.

Conversely, if $\tilde{X}\in \mathfrak X _L(G)$ and $\tilde{\mu}$
is the 1-form on $G$ defined by $\tilde{\mu}=-i(\tilde{X})\Omega$
then, from (\ref{lagrangiano}), it follows that $e^{\sigma}\s
^\ast \mu=\tilde{\mu}$, $\mu$ being the 1-form on $M$ given by
$\mu =\epsilon ^\ast \tilde{\mu}$. This implies that $\Psi (\mu
)=\tilde{X}$.

On the other hand, using that the map $\sostP :\Omega ^1(G)\to
\mathfrak X (G)$ is an isomorphism of $C^\infty (G,\R)$-modules,
we conclude that $\Psi$ is an isomorphism of $C^\infty
(M,\R)$-modules.

Now, suppose that $X,Y\in \Gamma (AG)$. We have that the
left-invariant vector field $\izq{X}$ is $\s$-projectable to a
vector field $a(X)$ on $M$. In addition, if $\mu$ and $\nu$ are
the 1-forms on $M$ satisfying $\Psi (\mu )=\izq{X}$ and $\Psi (\nu
)=\izq{Y}$, then a long computation, using (\ref{0'}) and
Definition \ref{def-lcs}, shows that
\[
a(X)= \widetilde{\#}_{(\Lambda _0,E_0,\theta _0)}(\mu ),\qquad
i([\izq{X},\izq{Y}])\Omega =-e^{\sigma }\s ^\ast \lcf \mu ,\nu
\rcf _{(\Lambda _0 ,E_0,\theta _0)}.
\]
This ends the proof of our result.\QED
\begin{remark}
{\rm  Let $(\gr ,\Omega )$ be a symplectic groupoid. Then, the
Jacobi structure on $M$ is Poisson, that is, $E_0=0$ (see Section
5.2 in \cite{IM}) and the 1-form $\theta _0$ on $M$ identically
vanishes. Thus (see Theorem \ref{algebroide}), $AG$ is isomorphic
to the cotangent Lie algebroid $T^\ast M$. This result was proved
in \cite{CDW}. }
\end{remark}

\section*{Acknowledgments}
Research partially supported by DGICYT grant BFM 2000-0808. D.
Iglesias-Ponte wishes to thank the Spanish Ministry of Education
and Culture for an FPU grant.

%
%

\label{endxyzt}
\end{document}